\documentclass[11pt]{article}

\usepackage{amsmath,amssymb,amstext,amsthm}
\usepackage{geometry}
\usepackage{fancyhdr}
\usepackage[font=small,labelfont=bf, width=.85\textwidth]{caption}
\usepackage[pdftex]{graphicx}
\usepackage{xcolor}

\usepackage[
pdfauthor={},
pdftitle={},
pdfstartview=XYZ,
bookmarks=true,
colorlinks=true,
linkcolor=blue,
urlcolor=blue,
citecolor=blue,
bookmarks=true,
linktocpage=true,
hyperindex=true
]{hyperref}

\newtheorem{theorem}{Theorem}

\newtheorem{clm}{Claim}

\newenvironment{proofc}{{\noindent \textit{Proof of Claim.}}}
{\hfill $\Box$}

\begin{document}

\title{Cliques and High Odd Holes in Graphs with Chromatic Number Equal to Maximum Degree}
\author{
Rachel Galindo\footnote{Auburn University, Department of Mathematics and Statistics, Auburn U.S.A.
  Email: {\tt rcg0036@auburn.edu}.}
\qquad
Jessica McDonald\footnote{Auburn University, Department of Mathematics and Statistics, Auburn U.S.A.
  Email: {\tt mcdonald@auburn.edu}.
	Supported in part by Simons Foundation Grant \#845698  }
\qquad
Songling Shan \footnote{Auburn University, Department of Mathematics and Statistics, Auburn U.S.A
Email: {\tt szs0398@auburn.edu}. Supported in part by NSF grant
			DMS-2345869.}}

\date{}

\maketitle

\begin{abstract}
We give a uniform and self-contained proof that if $G$ is a connected graph with $\chi(G) = \Delta(G)$ and $G\neq \overline{C_7}$, then $G$ contains either $K_{\Delta(G)}$ or an odd hole where every vertex has degree at least $\Delta(G)-1$ in $G$. This was previously proved in series of two papers by Chen, Lan, Lin, and Zhou, who used the Strong Perfect Graph Theorem for the cases $\Delta(G)=4, 5, 6$.
\end{abstract}

\vspace*{.3in}

\section{Introduction}

In this paper all graphs are simple; we follow \cite{West} for standard terms not defined here.

It is easy to show that every graph $G$ satisfies $\chi(G)\leq \Delta(G)+1$, where  $\chi(G), \Delta(G)$ denote the chromatic number and maximum degree of $G$, respectively. Brooks' Theorem \cite{Brooks} says that if a connected graph $G$ has $\chi(G)=\Delta(G)+1$, then $G$ is either  $K_{\Delta(G)+1}$ or an odd cycle. In 2023, Chen, Lan, Lin, and Zhou \cite{ChenLanLinZhou} proved that if $\chi(G)=\Delta(G)\geq 7$, then $G$ contains either  $K_{\Delta(G)}$ or an odd hole, that is, a chordless odd cycle of length at least five. Although not stated in the paper, their proof actually finds an odd hole where all vertices have degree at least $\Delta(G)-1$ in $G$. We term such a subgraph a \emph{high odd hole}.

\begin{theorem} [Chen, Lan, Lin, Zhou \cite{ChenLanLinZhou}]\label{thm: ChenLanLinZhou}
    Let $G$ be a graph with $\chi(G) = \Delta(G) \geq 7$. Then $G$ contains either a $K_{\Delta(G)}$ or a high odd hole.
\end{theorem}

It is easy to show that Theorem \ref{thm: ChenLanLinZhou} also holds for $\Delta(G) \leq 3$. (A triangle-free graph $G$ with $\chi(G) = 3$ must contain an odd cycle of length at least 5, the shortest of which is an odd hole.) On the other hand, the graph $\overline{C_7}$ (the complement of the cycle on seven vertices) has chromatic number and maximum degree both equal to 4, but contains no odd holes and no copies of $K_4$; this was first pointed out to us by Xie \cite{Personalcommunicationxie}. In fact, it turns out that $\overline{C_7}$ is the unique exceptional graph.

\begin{theorem} \label{thm: main}
Let $G$ be a connected graph with $\chi(G) = \Delta(G)$ and $G\neq \overline{C_7}$. Then $G$ contains either a $K_{\Delta(G)}$ or a high odd hole.
\end{theorem}

Theorem \ref{thm: main} is the focus of this paper. In fact, Theorem \ref{thm: main} was very recently proved by Chen, Lan, Lin, and Zhou \cite{CHEN2025}; our proof is independent and quite different for the cases $\Delta(G)=4, 5, 6$. Chen, Lan, Lin and Zhou use the Strong Perfect Graph Theorem for these smaller cases. Our proof of Theorem \ref{thm: main} is uniform and self-contained.

The Borodin-Kostochka Conjecture \cite{BK} from 1977 posits that if $G$ is a graph with $\chi(G)=\Delta(G)\geq 9$, then $G$ contains a $K_{\Delta(G)}$. Despite a large literature of results, the conjecture remains open; the reader is referred to \cite{CranRab,Jensen} for details about its history and currently known partial results. Here we note that when $\Delta(G)\geq 9$, Theorems \ref{thm: ChenLanLinZhou} and \ref{thm: main} are a (significant) weakening of the Borodin-Kostochka Conjecture.  When $3\leq \Delta(G)\leq 8$ there are well-known examples showing that $G$ need not contain a $K_{\Delta(G)}$ when $\chi(G)=\Delta(G)$.  In particular, there are some examples where every odd hole has at least two vertices of degree $\Delta(G)-1$ (eg. by Catlin \cite{Catlin} for $\Delta(G)=6,7$ and by Kostochka, Rabern, and Steibitz \cite{DefsPaper} for $\Delta(G)=5$). Hence we see that the odd hole is necessary in the statement of Theorems \ref{thm: ChenLanLinZhou} and \ref{thm: main}, and that its degree bound is as high as we can hope for.

\section{Proof of Theorem \ref{thm: main}}

Let $G$ be a connected graph with $\chi(G) = \Delta(G) \geq 4$; let $\Delta=\Delta(G)$. We suppose that $G$ contains no $K_{\Delta}$ and no high odd hole, and we will show that $G=\overline{C_7}$. Let $H$ be a vertex-critical subgraph of $G$ with $\chi(H)=\chi(G)$, obtained by deleting vertices from $G$. Note that $H$ is an induced subgraph of $G$, and moreover that $\delta(H)\geq \Delta(H)-1$ since it is vertex-critical. If $\Delta(H) = \Delta - 1,$ then $\chi(H)=\Delta(H)+1$, so Brooks' Theorem implies that  $H=K_{\Delta(H)+1}=K_{\Delta}$ (since $\Delta(H)\geq 3$), contradiction. So we know that $\Delta(H)=\Delta$.  We prove the following five claims.

\begin{clm}\label{clm:adj0} Let $x, y, z\in V(H)$ and let $\varphi$ be a $(\Delta-1)$-coloring of $H-x$. If $x\sim y, z$ and $\varphi(y), \varphi(z)$ are distinct colors that are not used on any other neighbors of $x$, then  $y\sim z$.
\end{clm}

\begin{proofc} Consider the maximal $(\varphi(y), \varphi(z))$-alternating subgraph $H'\subseteq H$ starting at $y$.
Note that $H'$ must contain $z$, since if not by swapping the two colors on all of $H'$ we get a new $(\Delta-1)$-coloring of $G-x$ in which color $\varphi(y)$ is missing from the neighborhood of $x$, contradicting $\chi(H)=\Delta$. Let $P$ be the shortest path from $y$ to $z$ in $H'$. If $y\not \sim z$ then $P$ along with $x$ is an odd hole in $H$, and also in $G$ since $H$ is induced. Since every vertex in $H$ has degree at least $\Delta(H)-1=\Delta-1$, this is a high odd hole in $G$, contradiction.
\end{proofc}

\begin{clm}\label{clm:reg} We must have $G=H$ and this graph is $\Delta$-regular.
\end{clm}
\begin{proofc} Let $v\in V(H)$ and let $\varphi$ be a $(\Delta-1)$-coloring of $H-v$. Since $\varphi$ cannot be extended to $v$, each of the colors $1, 2, \ldots, \Delta-1$ must occur on the neighbors of $v$; one color may appear twice if $d_H(v)=\Delta$. If $d_H(v)=\Delta-1$ then by Claim \ref{clm:adj0} the neighbors of $v$ induce a $K_{\Delta-1}$ in $H$, and adding $v$ gives a $K_{\Delta}$, contradiction. So $H$ must be $\Delta$-regular. Since $G$ is connected, this implies $H=G$.
\end{proofc}

\begin{clm}\label{clm:nbhd}
For any vertex $v\in V(G)$, its neighborhood $N(v)$ can be partitioned into $A_v, B_v$ with:
\begin{itemize}
\item[(a)] $A_v$ is an independent set of size 2;
\item[(b)] $B_v$ is a clique of size $\Delta-2$, and;
\item[(c)] every vertex in $B_v$ is adjacent to at least one vertex in $A_v$.
\end{itemize}
\end{clm}

\begin{proofc}
Let $v\in V(G)$ and let $\varphi$ be a $(\Delta-1)$-coloring of $G-v$. We may assume without loss that the neighbors of $v$ are $v_0, v_1, \ldots v_{\Delta-1}$ with $\varphi(v_i)=i$ for $1\leq i\leq \Delta-1$ and $\varphi(v_0)=1$, by Claim \ref{clm:reg} and by the fact that $\varphi$ cannot be extended to $v$. Then $v_0\not\sim v_1$ and Claim \ref{clm:adj0} tells us that the vertices $v_2, \ldots, v_{\Delta-1}$ induce a $K_{\Delta-2}$. Let $A_v=\{v_0, v_1\}$ and $B_v=\{v_2, \ldots, v_{\Delta-1}\}$. It then remains only to prove (c). To this end, we fix $i\in\{2, \ldots, {\Delta-1}\}$ and show that $v_i$ is adjacent to at least one of $v_0, v_1$.

Consider the maximal $(1,i)$-alternating subgraph $G_{1i}$ starting at $v_i$. If $G_{1i}$ contains neither of $v_0, v_1$, then by swapping the two colours on all of $G_{1i}$ we get a new $(\Delta-1)$-coloring of $G-v$ in which colour $i$ is missing from the neighborhood of $v$, contradicting $\chi(H)=\Delta$. So $G_{1i}$ must contain at least one of $v_0, v_1$.  Let $P$ be a shortest path  in $G_{1i}$ from $v_i$ to $\{v_0, v_1\}$. If $P$ is a single edge, then $v_i$ is adjacent to at least one of $v_0, v_1$ as desired. Otherwise, $P$ along with $v$ is an odd hole in $G$, contradiction.
\end{proofc}

\begin{clm}\label{clm:NaBv}
Let $v\in V(G)$ and let $a \in A_v$. Then $|N(a) \cap  B_v|=\Delta-3$.
\end{clm}

\begin{proofc} If $|N(a)\cap B_v|=|B_v|=\Delta-2$, then $B_v\cup\{a, v\}$ induces a $K_{\Delta}$ in $G$, contradiction. On the other hand, if $|N(a)\cap B_v|=0$, then by Claim \ref{clm:nbhd} there exists  $a' \in A_v\setminus\{a\}$ with $|N(a')\cap B_v|=|B_v|=\Delta-2$, contradiction. So we get that
$$1\leq |N(a) \cap  B_v| \leq \Delta-3.$$
\indent Suppose first that there exists $b\in B_v\cap B_{a}$ (where $B_a$ is obtained by applying Claim \ref{clm:nbhd} to the vertex $a$). Since $b\in B_v$ we know that $b$ has $\Delta-3$ neighbors in $B_v$, plus it is adjacent to both $v$ and to $a$ (since $b\in B_a$). This accounts for $\Delta-1$ out of its $\Delta$ neighbors. Since $b\in B_{a}$, this means that $B_{a}$ contains at most one vertex outside of $B_v\cup \{v\}$. If $B_a$ contains a vertex outside of $B_v\cup \{v\}$, then $v\not\in B_a$ (if $v \in B_a$ then $v$ has $|N(v)| > \Delta$, contradiction), so $a$ has at least $\Delta-3$ neighbors in $B_v$, as desired. On the other hand, if $B_a$ contains no vertices outside of $B_v\cup \{v\}$, then again $a$ must have at least $\Delta-3$ neighbors in $B_v$.

We may now assume that $B_v\cap B_{a}=\emptyset$. Since $|N(a)\setminus B_a|=2$, and these two vertices are non-adjacent by Claim \ref{clm:nbhd}, this implies that $|N(a) \cap  B_v|=1$, say $w\in N(a)\cap B_v$. Then $w\in A_a$. Since $w\sim v$ this means that $v\not\in A_a$. But then since $a\sim v$ we must have $v\in B_a$. But then the $\Delta-2\geq 2$ vertices of $B_a$ must all be in  $N(v)\cup \{v\}$, which is impossible since $B_v \cap  B_a=\emptyset$ and $A_v\setminus \{a\}$ contains only one vertex, which is not adjacent to $a$ (by Claim \ref{clm:nbhd}), contradiction.
\end{proofc}

\begin{figure}
    \centering
    \includegraphics[width=.65\linewidth]{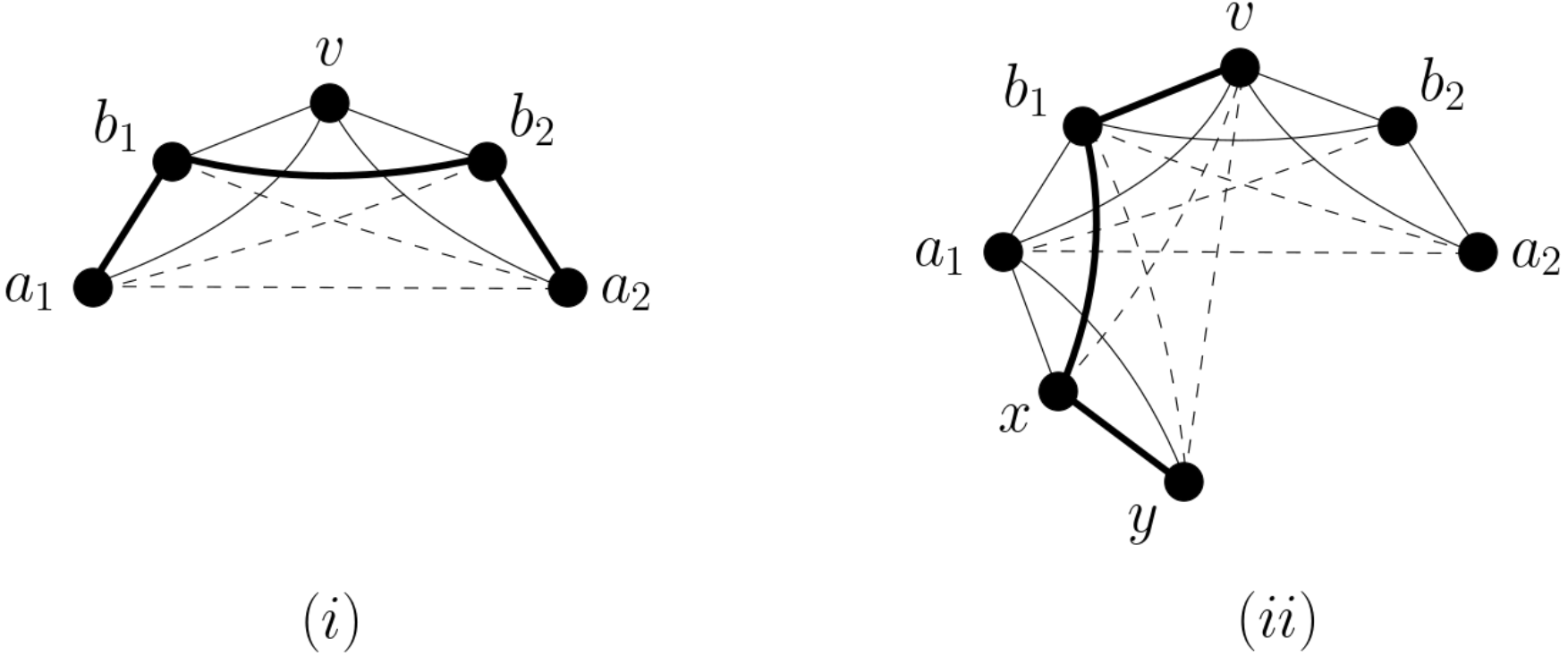}
    \caption{Every vertex $v\in V(G)$ has four neighbors that look like image $(i)$, as described in Claim \ref{clm:4nbrs}: the path $(a_1, b_1, b_2, a_2)$ is indicated in bold, and the dotted lines indicate non-adjacency. When $\Delta=4$, we may apply Claim \ref{clm:4nbrs} to both $v$ and $a_1$ to get the structure in (ii).}
    \label{fig:Fig1}
\end{figure}

\begin{clm}\label{clm:4nbrs}
Let $v\in V(G)$. Then there exist four distinct vertices $a_1, a_2, b_1, b_2\in N(v)$ such that $G[\{a_1, a_2, b_1, b_2\}]$ is precisely the path $(a_1,b_1, b_2, a_2)$. See Figure \ref{fig:Fig1}(i).
\end{clm}

\begin{proofc}
Consider the the partition of  $A_v, B_v$ guaranteed by Claim \ref{clm:nbhd}, and let $A_v=\{a_1, a_2\}$. By Claim \ref{clm:NaBv} there is a unique vertex in $B_v$ that is non-adjacent to $a_1$, say $b_2$, and a unique vertex in $B_v$ that is non-adjacent to $a_2$, say $b_1$. We know that $b_1\neq b_2$ by Claim \ref{clm:nbhd}(c). Since $A_v$ is an independent set, and $B_v$ is a clique, we have now completely determined the edges in $G[\{a_1, a_2, b_1, b_2\}]$, as described.
\end{proofc}

\vspace*{.1in}

Fix $v\in V(G)$, and consider the vertices $a_1, a_2, b_1, b_2$ guaranteed by Claim \ref{clm:4nbrs}.

Suppose first that $\Delta\geq 5$. Then there exists $z\in N(v)\setminus \{a_1, a_2, b_1, b_2\}$. By the proof of Claim \ref{clm:NaBv} we may assume that $A_v=\{a_1, a_2\}$ and $b_1, b_2, z\in B_v$. Note that $z$ has $\Delta-3$ neighbors in $B_v$, plus it is adjacent to all of $v, a_1, a_2$ by Claims \ref{clm:NaBv} and \ref{clm:4nbrs}. This accounts for all $\Delta$ of its neighbors. On the other hand, $a_1$ has exactly $\Delta-2$ neighbors in $N(v)\cup\{v\}$ and two from outside this set, say $s,t$.

If $z\in B_{a_1}$ then since $z\not\sim s, t$, we get that $s, t\not\in B_{a_1}$. So we must have $A_{a_1}=\{s, t\}$. But then by Claim \ref{clm:nbhd}(c) the vertex $z$ must be adjacent to at least one of $s, t$, contradiction.

Suppose now that $z\in A_{a_1}$. Since $|A_{a_1}|=2$, we may assume without loss of generality that $t\in B_{a_1}$. Since $t\not\sim v$ this means that $v\not\in B_{a_1}$ and hence that $A_{a_1}=\{z, v\}$. But since $v\sim z$, this contradicts Claim \ref{clm:nbhd}(a).

We may now assume that $\Delta=4$. In particular, this means that $N(v)=\{a_1, a_2, b_1, b_2\}$ and by Claim \ref{clm:4nbrs} $N(v)\cup \{v\}$ induces precisely the graph depicted in Figure 1(i). However, we may also apply Claim \ref{clm:4nbrs} to the vertex $a_1$ (in place of $v$).  When we do this, we get that there exists vertices $x, y$ such that $N(a_1)=\{v, b_1, x, y\}$ and $G[N(a_1)]$ is precisely the path $(v, b_1, x, y)$. We know that $x, y\not\in\{a_2, b_2\}$ since $v\sim a_2, b_2$ while $v\not\sim x, y$.  See Figure \ref{fig:Fig1}(ii).

\begin{figure}
    \centering
    \includegraphics[width=0.75\linewidth]{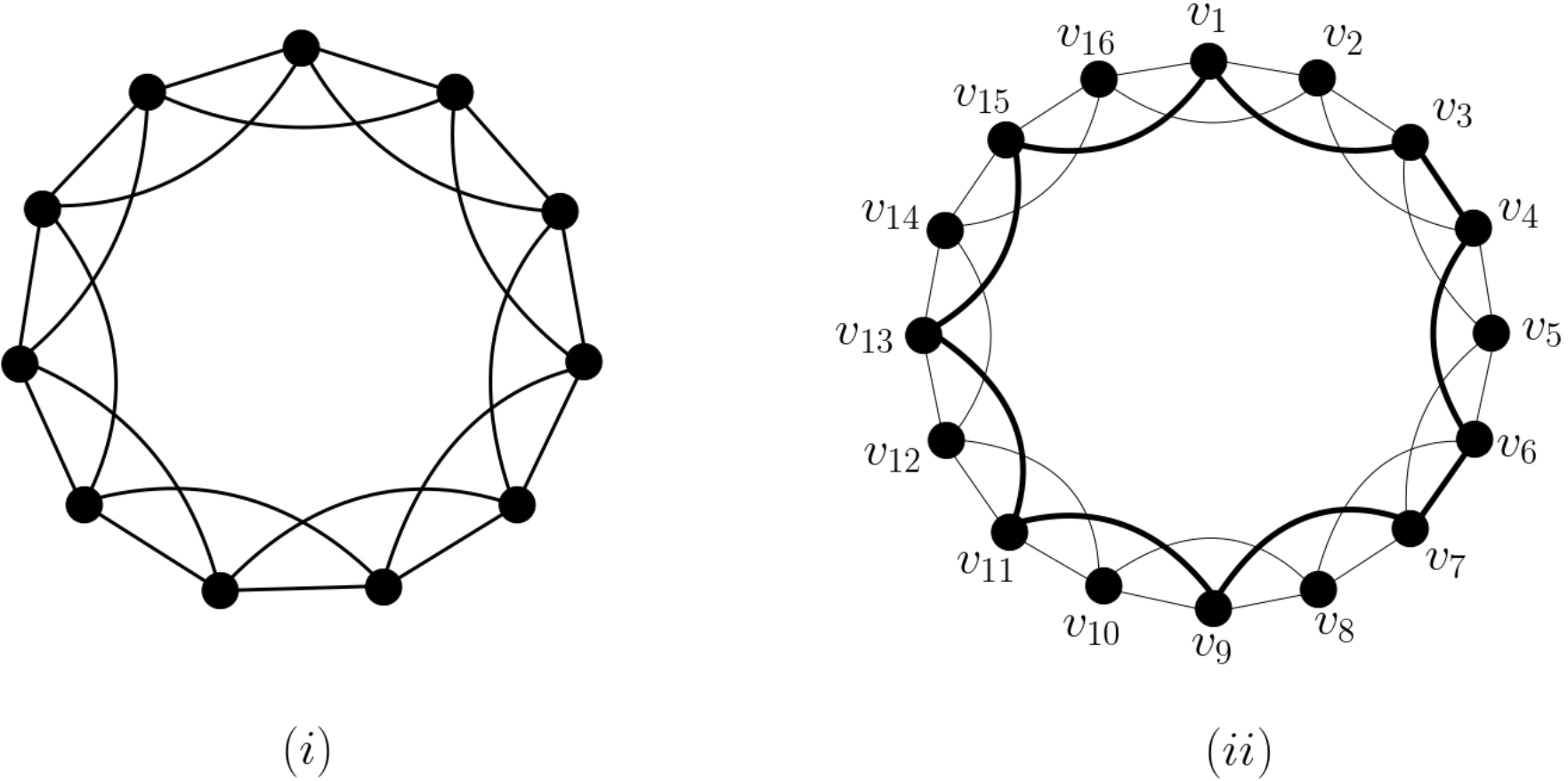}
    \caption{(i) The graph obtained from the cycle $C_{11}$ by adding an edge between all its distance-2 vertices. (ii) In the case $n=16$ we get the odd-hole $L$ indicated with bold edges.}
    \label{fig:Fig2}
\end{figure}

Note that at this point, we have shown that $G$ contains at least 7 vertices, and the neigborhoods of $v, b_1, a_1$ have been completely determined-- they will not be adjacent to any additional vertices. We can now continue in this process and apply Claim \ref{clm:4nbrs} to $y$, yielding vertices $u, v$ such that $N(y)=\{a_1, x, u, w\}$ and $G[N(y)]$ is precisely the path $(a_1,x, u, w)$. Now however, there are two possibilities, considering that Claim \ref{clm:4nbrs} can be applied also to $b_2, a_2$: either $u, w\not\in \{a_2, b_2\}$, or $u=a_2$ and $w=b_2$.  In the former case we continue our argument by applying Claim \ref{clm:4nbrs} to $w$ to get two new vertices, and so on. Eventually however, since the graph is 4-regular, we will get that the two new vertices are indeed $a_2, b_2$. Hence we get that the graph $G$ is be obtained from a cycle $C_n$ with $n\geq 7$ by adding an edge between all its distance-2 vertices (we call these added edges the \emph{distance-2} edges of $G$, and we call the others the \emph{distance-1} edges of $G$). See Figure \ref{fig:Fig2}(i).

If $n\equiv 0$ (mod 3), then $G$ has a 3-coloring obtained by using the colors 1, 2, 3 in sequence around $C_n$, contradiction.

Suppose now that that $n \equiv 2$ (mod 3). Choose any vertex $v \in V(G)$ and let $G' = G - v$. Then $\chi(G') = 3$ since $G$ is critical, so it is possible to assign a 3-coloring to $G'$. To this end, choose a vertex $w$ in $G'$ which was adjacent to $v$ in $G$ and assign it color 1. As $w$ is in a $K_3$ with the next two consecutive vertices on $C_n$, assign these vertices colors 2 and 3, respectively, noting that the coloring is without loss of generality at this point.  Moreover, the structure of $G'$ forces us to repeat the pattern 1, 2, 3 as we go around the cycle. Suppose that $x$ is the last vertex  to receive a color as we travel around the cycle this way. Then since $n \equiv 2$ (mod 3), $G'$ has $n-1\equiv 1$ (mod 3) vertices, and we know that $x$ must receive color 1. However we also have $x\sim w$ which was also assigned color 1, contradicting the fact that $\chi(G')=3$.

We may now assume that $n \equiv 1$ (mod 3); let $n = 3k + 1$. If $k=2$ then $G=\overline{C_7}$, as desired, so we may assume that $k\geq 3$. Label the vertices of $G$ as $v_1, v_2, \dots v_{3k + 1}$, moving clockwise around the cycle $C_n$. We claim that $G$ contains an odd hole. We build this odd hole $L$ by beginning at $v_1$ and taking the distance-2 edge to $v_3$. We then continue to move clockwise around the cycle, alternating taking a distance-1 edge and a distance-2 edge, until we have taken a total of $k-3$ distance-2 edges and $k-3$ distance-1 edges. Since we started at $v_1$, and $1+2(k-3)+(k-3)=3k-8$, this means we have stopped at the vertex $v_{3k-8}$. We complete the cycle $L$ by taking five distance-2 edges back to $v_{3k+1}$ (noting that $3k-8 +2(5)=(3k+1)+1$). See Figure \ref{fig:Fig2}(ii). The cycle $L$ is a hole since we did not take any consecutive distance-1 edges. Moreover, $L$ has total length $(k-2)+(k-3)+4=2k-1$, which is odd. So we have found an odd hole in $G$, contradiction.\qed

\bibliographystyle{amsplain}
\bibliography{UpdatedCliquesAndHighOddHoles}

\end{document}